\documentclass[graybox]{svmult}

\usepackage{mathptmx}       
\usepackage{helvet}         
\usepackage{courier}        
\usepackage{makeidx}         
\usepackage{graphicx}        
\usepackage{multicol}        

\usepackage{latexsym,amsmath,amsfonts}
\usepackage{enumerate}

\usepackage{epsfig}
\usepackage{verbatim}

\makeindex             

\begin{document}

\title*{Multilevel Monte Carlo methods}
\author{Michael B.~Giles}
\institute{Michael B.~Giles \at Mathematical Institute,
                  University of Oxford,\ \email{mike.giles@maths.ox.ac.uk}}

%
%
\maketitle

\abstract{
The author's presentation of multilevel Monte Carlo path simulation
at the MCQMC 2006 conference stimulated a lot of research into 
multilevel Monte Carlo methods.  This paper reviews the progress 
since then, emphasising the simplicity, flexibility and generality of 
the multilevel Monte Carlo approach.  It also offers a few original 
ideas and suggests areas for future research.
}

\section{Introduction}

\subsection{Control variates and two-level MLMC}

One of the classic approaches to Monte Carlo variance reduction is 
through the use of a control variate.
Suppose we wish to estimate $\mathbb{E}[f]$, and there is a control variate 
$g$ which is well correlated to $f$ and has a known expectation $\mathbb{E}[g]$.
In that case, we can use the following unbiased estimator for
$\mathbb{E}[f]$:
\[
N^{-1} \sum_{n=1}^N
\left\{ f^{(n)} - \lambda \left( g^{(n)} - \mathbb{E}[g] \right) \right\}.
\]
The optimal value for $\lambda$ is $\rho \sqrt{\mathbb{V}[f]\, /\, \mathbb{V}[g]}$,
where $\rho$ is the correlation between $f$ and $g$, and the variance 
of the control variate estimator is reduced by factor $1\!-\!\rho^2$ 
compared to the standard estimator.

A two-level version of MLMC (multilevel Monte Carlo) is very similar.
If we want to estimate $\mathbb{E}[P_1]$ but it is much cheaper to
simulate $P_0 \approx P_1$, then since
\[
\mathbb{E}[P_1] = \mathbb{E}[P_0] + \mathbb{E}[P_1-P_0]
\]
we can use the unbiased two-level estimator
\[
N_0^{-1} \sum_{n=1}^{N_0} P_0^{(n)} \ + \ 
N_1^{-1} \sum_{n=1}^{N_1} \left(P_1^{(n)} - P_0^{(n)}\right).
\]
Here $P_1^{(n)} \!-\! P_0^{(n)}$ represents the difference between 
$P_1$ and $P_0$ for the same underlying stochastic sample, so that
$P_1^{(n)} \!-\! P_0^{(n)}$ is small and has a small variance; the
precise construction depends on the application and various examples
will be shown later.
The two key differences from the control variate approach are that 
the value of $\mathbb{E}[P_0]$ is not known, so has to be estimated, and 
we use $\lambda=1$.

If we define $C_0$ and $C_1$ to be the cost of computing a single 
sample of $P_0$ and $P_1\!-\!P_0$, respectively, 
then the total cost is
$
N_0\, C_0 \!+\! N_1\, C_1,
$
and if $V_0$ and $V_1$ are the variance of $P_0$ and $P_1\!-\!P_0$, 
then the overall variance is
$
N_0^{-1} V_0 + N_1^{-1} V_1,
$
assuming that $\displaystyle \sum_{n=1}^{N_0} P_0^{(n)}$ and
$\displaystyle \sum_{n=1}^{N_1} \left(P_1^{(n)} - P_0^{(n)}\right)$
use independent samples.

Hence, treating the integers $N_0, N_1$ as real variables and
performing a constrained minimisation using a Lagrange multiplier, 
the variance is minimised for a fixed cost by choosing
$
N_1 /\, N_0 = \sqrt{V_1/C_1} \ /  \sqrt{V_0/C_0}.
$

\subsection{Multilevel Monte Carlo}

\label{sec:unbiased_MLMC}

The full multilevel generalisation is quite natural: given 
a sequence $P_0, P_1, \ldots,$ which approximates $P_L$ with 
increasing accuracy, but also increasing cost, we have the 
simple identity
\[
\mathbb{E}[P_L] = \mathbb{E}[P_0] + \sum_{\ell=1}^L \mathbb{E}[P_\ell-P_{\ell-1}],
\]
and therefore we can use the following unbiased estimator
for $\mathbb{E}[P_L]$,
\[
N_0^{-1} \sum_{n=1}^{N_0} P_0^{(0,n)} \ + \ 
 \sum_{\ell=1}^L \left\{
N_\ell^{-1} \sum_{n=1}^{N_\ell} \left(P_\ell^{(\ell,n)} - P_{\ell-1}^{(\ell,n)}\right)
\right\}
\]
with the inclusion of the level $\ell$ in the superscript $(\ell,n)$ 
indicating that the samples used at each level of correction are independent.

If we define $C_0, V_0$ to be the cost and variance of one sample 
of $P_0$, and $C_\ell, V_\ell$ to be the cost and variance of one 
sample of $P_\ell\!-\!P_{\ell-1}$, then the overall cost and variance
of the multilevel estimator is
$\displaystyle
\sum_{\ell=0}^L N_\ell\ C_\ell\ 
$
and 
$\ \displaystyle
\sum_{\ell=0}^L N_\ell^{-1} \ V_\ell
$,
respectively.

For a fixed cost, the variance is minimised by choosing
$
N_\ell = \lambda \sqrt{V_\ell \, / \, C_\ell}
$
for some value of the Lagrange multiplier $\lambda$.  
In particular, to achieve an overall variance of $\varepsilon^2$ 
requires that
$
\lambda = \varepsilon^{-2}\ \sum_{\ell=0}^L \sqrt{V_\ell \ C_\ell}.
$
The total computational cost is then
\begin{equation}
C = 
\varepsilon^{-2} \left(\sum_{\ell=0}^L \sqrt{V_\ell \ C_\ell} \right)^2.
\label{eq:total_cost}
\end{equation}

It is important to note whether the product $V_\ell \ C_\ell$ increases 
or decreases with $\ell$, i.e.~whether or not the cost increases with
level faster than the variance decreases.  If it increases with level,
so that the dominant contribution to the cost comes from $V_L \, C_L$
then we have $C \approx \varepsilon^{-2} V_L \, C_L$,
whereas if it decreases and the dominant contribution comes from  
$V_0 \, C_0$ then $C \approx \varepsilon^{-2} V_0 \, C_0$.
This contrasts to the standard MC cost of approximately
$\varepsilon^{-2} V_0 \, C_L$, assuming that the cost of computing $P_L$ 
is similar to the cost of computing $P_L \!-\!P_{L-1}$, and that
$\mathbb{V}[P_L] \approx \mathbb{V}[P_0]$.  This shows that in the first case the
MLMC cost is reduced by factor $V_L/V_0$, corresponding to the 
ratio of the variances $\mathbb{V}[P_L \!-\!P_{L-1}]$ and $\mathbb{V}[P_L]$,
whereas in the second case it is reduced by factor $C_0/C_L$, 
the ratio of the costs of computing $P_0$ and $P_L\!-\!P_{L-1}$.
If the product $V_\ell \ C_\ell$ does not vary with level, then the 
total cost is $\varepsilon^{-2}L^2\,V_0 \, C_0 =  \varepsilon^{-2}L^2\,V_L \, C_L$.

\subsection{Earlier related work}

Prior to the author's first publications \cite{giles08b,giles08}
on MLMC for Brownian path simulations, Heinrich developed a 
multilevel Monte Carlo method for parametric integration, 
the evaluation of functionals arising from the solution of 
integral equations, and weakly singular integral operators
\cite{heinrich98,heinrich00,heinrich01,heinrich06,hs99}.  
Parametric integration concerns the estimation of $\mathbb{E}[f(x,\lambda)]$ 
where $x$ is a finite-dimensional random variable and $\lambda$ is a 
parameter.  In the simplest case in which $\lambda$ is a real variable 
in the range $[0,1]$, having estimated the value of $\mathbb{E}[f(x,0)]$ and 
$\mathbb{E}[f(x,1)]$, one can use $\frac{1}{2}(f(x,0)+f(x,1))$ as a control 
variate when estimating the value of $\mathbb{E}[f(x,\frac{1}{2})]$.
This approach can then be applied recursively for other intermediate
values of $\lambda$, yielding large savings if $f(x,\lambda)$ is 
sufficiently smooth with respect to $\lambda$.  Although this
does not quite fit into the general MLMC form given in the previous 
section, the recursive control variate approach is very similar and 
the complexity analysis is also very similar to the analysis to be 
presented in the next section.

Although not so clearly related, there are papers by Brandt 
{\it et al} \cite{bgr94,bi03} which combine Monte Carlo techniques
with multigrid ideas in determining thermodynamic limits in 
statistical physics applications.  It is the multigrid ideas of 
Brandt and others for the iterative solution of systems of equations 
which were the inspiration for the author in developing the MLMC method
for SDE path simulation.

In 2005, Kebaier \cite{kebaier05} developed a two-level approach
for path simulation which is very similar to the author's approach
presented in the next section.  The only differences are the use 
of only two levels, and the use of a general multiplicative factor 
as in the standard control variate approach.  A similar multilevel 
approach was under development at the same time by Speight, but was 
not published until later \cite{speight09,speight10}.


\section{MLMC theorem}

In the Introduction, we considered the case of a general multilevel
method in which the output $P_L$ on the finest level corresponds 
to the quantity of interest. However, in many infinite-dimensional
applications, such as in SDEs and SPDEs, the output $P_\ell$ on 
level $\ell$ is an approximation to a random variable $P$.  
In this case, the mean square error (MSE) has the usual decomposition 
into the total variance of the multilevel estimator, plus the square 
of the bias $(\mathbb{E}[P_L \!-\!P])^2 $.  
To achieve an MSE which is less than $\varepsilon^2$, it is sufficient to 
ensure that each of these terms is less than
${\textstyle \frac{1}{2}}\varepsilon^2$.
This leads to the following theorem:

\begin{theorem}
\label{thm:MLMC}
Let $P$ denote a random variable, and let $P_\ell$ denote the 
corresponding level $\ell$ numerical approximation.

If there exist independent estimators $Y_\ell$ based on $N_\ell$ 
Monte Carlo samples, and positive constants 
$\alpha, \beta, \gamma, c_1, c_2, c_3$ such that 
$\alpha\!\geq\!{\textstyle \frac{1}{2}}\,\min(\beta,\gamma)$ and
\begin{itemize}
\item[i)] ~
$\displaystyle
\left| \mathbb{E}[P_\ell \!-\! P] \right|\ \leq\ c_1\, 2^{-\alpha\, \ell}
$
\item[ii)] ~
$\displaystyle
\mathbb{E}[Y_\ell]\ = \left\{ \begin{array}{ll}
\mathbb{E}[P_0],                     &~~ \ell=0 \\[0.1in]
\mathbb{E}[P_\ell \!-\! P_{\ell-1}], &~~ \ell>0
\end{array}\right.
$
\item[iii)] ~
$\displaystyle
\mathbb{V}[Y_\ell]\ \leq\ c_2\, N_\ell^{-1}\, 2^{-\beta\, \ell}
$
\item[iv)] ~
$\displaystyle
\mathbb{E}[C_\ell]\ \leq\ c_3\, N_\ell\, 2^{\gamma\, \ell},
$
where $C_\ell$ is the computational complexity of $Y_\ell$
\end{itemize}
then there exists a positive constant $c_4$ such that for any 
$\varepsilon \!<\! e^{-1}$
there are values $L$ and $N_\ell$ for which the multilevel estimator
\[
Y = \sum_{\ell=0}^L Y_\ell,
\]
has a mean-square-error with bound
\[
MSE \equiv \mathbb{E}\left[ \left(Y - \mathbb{E}[P]\right)^2\right] < \varepsilon^2
\]
with a computational complexity $C$ with bound
\[
\mathbb{E}[C] \leq \left\{\begin{array}{ll}
c_4\, \varepsilon^{-2}              ,    & ~~ \beta>\gamma, \\[0.1in]
c_4\, \varepsilon^{-2} (\log \varepsilon)^2,    & ~~ \beta=\gamma, \\[0.1in]
c_4\, \varepsilon^{-2-(\gamma\!-\!\beta)/\alpha}, & ~~ \beta<\gamma.
\end{array}\right.
\]
\end{theorem}

The statement of the theorem is a slight generalisation of the 
original theorem in \cite{giles08}.  It corresponds to the theorem
and proof in \cite{cgst11}, except for the minor change to expected 
costs to allow for applications such as jump-diffusion modelling
in which the simulation cost of individual samples is itself random.

The theorem is based on the idea of a geometric progression in the 
levels of approximation, leading to the exponential decay in the 
weak error in condition {\it i)}, and the variance in condition 
{\it iii)}, as well as the exponential increase in the expected 
cost in condition {\it iv)}.  This geometric progression was based 
on experience with multigrid methods in the iterative solution of
large systems of linear equations, but it is worth noting that it 
is not necessarily the optimal choice in all circumstances.

The result of the theorem merits some discussion.  In the case 
$\beta>\gamma$, the dominant computational cost is on the coarsest 
levels where $C_\ell = O(1)$ and $O(\varepsilon^{-2})$ samples are 
required to achieve the desired accuracy.  This is the standard
result for a Monte Carlo approach using i.i.d.~samples; to do 
better would require an alternative approach such as the use
of Latin hypercube sampling or quasi-Monte Carlo methods.
In the case $\beta<\gamma$, the dominant computational cost is 
on the finest levels.  Because of condition {\it i)},
$2^{-\alpha L} = O(\varepsilon)$, and hence $C_L = O(\varepsilon^{-\gamma/\alpha})$.
If $\beta = 2 \alpha$, which is usually the largest possible value
for a given $\alpha$, for reasons explained below, 
then the total cost is $O(C_L)$ corresponding to $O(1)$ samples 
on the finest level, again the best that can be achieved.  
The dividing case $\beta=\gamma$ is the one 
for which both the computational effort, and the contributions to the 
overall variance, are spread approximately evenly across all of the 
levels; the $(\log\varepsilon)^2$ term corresponds to the $L^2$ factor in 
the corresponding discussion in section \ref{sec:unbiased_MLMC}.

The natural choice for the multilevel estimator is
\begin{equation}
\label{eq:natural}
Y_\ell = N_\ell^{-1}\ \sum_i\ 
P_\ell(\omega_i) \!-\! P_{\ell-1}(\omega_i),
\end{equation}
where $P_\ell(\omega_i)$ is the approximation to $P(\omega_i)$ 
on level $\ell$, and $P_{\ell-1}(\omega_i)$ is the corresponding 
approximation on level $\ell\!-\!1$ for the same underlying 
stochastic sample $\omega_i$.  Note that 
$\mathbb{V}[P_\ell \!-\!P_{\ell-1}]$ is usually similar in magnitude to
$\mathbb{E}[(P_\ell \!-\!P_{\ell-1})^2]$ which is greater than
$(\mathbb{E}[P_\ell \!-\!P_{\ell-1}])^2$; this implies that $\beta\leq 2 \alpha$
and hence the condition in the theorem that 
$\alpha \geq {\textstyle \frac{1}{2}}\,\min(\beta,\gamma)$ is satisfied.

However, the multilevel theorem allows for the use of other estimators,
provided they satisfy the restriction of condition {\it ii)} which
ensures that $\mathbb{E}[Y] = \mathbb{E}[P_L]$.
Two examples of this will be given later in the paper.  In the first,
slightly different numerical approximations are used for the coarse 
and fine paths in SDE simulations, giving
\[
Y_\ell = N_\ell^{-1}\ \sum_i\ 
P^f_\ell(\omega_i) \!-\! P^c_{\ell-1}(\omega_i).
\]
Provided $\mathbb{E}[P^f_\ell] = \mathbb{E}[P^c_\ell]$ so that the expectation 
on level $\ell$ is the same for the two approximations, then 
condition {\it ii)} is satisfied and no additional bias (other 
than the bias due to the approximation on the finest level) is 
introduced into the multilevel estimator.
The second example defines an antithetic $\omega_i^a$ with the 
same distribution as $\omega_i$, and then uses the multilevel
estimator 
\[
Y_\ell = N_\ell^{-1}\ \sum_i\ 
{\textstyle \frac{1}{2}} \left( P_\ell(\omega_i) \!+\! P_\ell(\omega^a_i) \right)
- P_{\ell-1}(\omega_i).
\]
Since $\mathbb{E}[P_\ell(\omega^a_i)] = \mathbb{E}[P_\ell(\omega_i)]$, then 
again condition {\it ii)} is satisfied.
In each case, the objective in constructing a more complex estimator 
is to achieve a greatly reduced variance $\mathbb{V}[Y_\ell]$ so that
fewer samples are required.

\section{SDEs}

\subsection{Euler discretisation}

The original multilevel path simulation paper \cite{giles08} 
treated SDEs using the simple Euler-Maruyama discretisation 
together with the natural multilevel estimator (\ref{eq:natural}).

Provided the SDE satisfies the usual conditions 
(see Theorem 10.2.2 in \cite{kp92}),
the strong error for the Euler discretisation with timestep $h$ 
is $O(h^{1/2})$, and therefore for Lipschitz payoff functions $P$ 
(such as European, Asian and lookback options in finance) 
the variance $V_\ell \equiv \mathbb{V}[P_\ell \!-\!P_{\ell-1}]$ is
$O(h_\ell)$.  If $h_\ell = 4^{-\ell} h_0$, as in \cite{giles08}, 
then this gives $\alpha\!=\!2$, $\beta\!=\!4$ and $\gamma\!=\!2$.
Alternatively, if $h_\ell = 2^{-\ell} h_0$, then
$\alpha\!=\!1$, $\beta\!=\!2$ and $\gamma\!=\!1$.
In either case, Theorem \ref{thm:MLMC} gives the complexity to 
achieve a root-mean-square error of $\varepsilon$ to be
$O(\varepsilon^{-2} (\log \varepsilon)^2)$, which is near-optimal as 
M{\"u}ller-Gronbach \& Ritter have proved an $O(\varepsilon^{-2})$
lower bound for the complexity \cite{mr09}.

For other payoff functions the complexity is higher. 
$V_\ell \approx O(h^{1/2})$ for the digital option which is a 
discontinuous function of the SDE solution at the final time, 
and the barrier option which depends discontinuously on the 
minimum or maximum value over the full time interval.  Loosely
speaking, this is because there is an $O(h^{1/2})$ probability of 
the coarse and fine paths being on opposite sides of the discontinuity, 
and in such cases there is an $O(1)$ difference in the payoff.
Currently, there is no known ``fix'' for this for the Euler-Maruyama
discretisation; we will return to this issue for the Milstein 
discretisation when there are ways of improving the situation.

Table \ref{tab:EM-Milstein} summarises the observed variance convergence 
rate in numerical experiments for the different options, and the 
theoretical results which have been obtained; the digital option
analysis is due to Avikainen \cite{avikainen09} while the others 
are due to Giles, Higham \& Mao \cite{ghm09}.  Although the analysis
in some of these cases is for one-dimensional SDEs, it also applies
to multi-dimensional SDEs \cite{giles09}.

\begin{table}[b!]
\begin{center}
\begin{tabular}{|l|l|l|l|l|}
\hline & \multicolumn{2}{c|}{Euler-Maruyama} &  \multicolumn{2}{c|}{Milstein} \\ 
option & numerics   & analysis & numerics  & analysis  \\ \hline
Lipschitz & $O(h)$       & $O(h)$              & $O(h^2)$ & $O(h^2)$ \\
Asian     & $O(h)$       & $O(h)$              & $O(h^2)$ & $O(h^2)$ \\
lookback  & $O(h)$       & $O(h)$              & $O(h^2)$ & $o(h^{2-\delta})$ \\
barrier   & $O(h^{1/2})$ & $o(h^{1/2-\delta})$ & $O(h^{3/2})$ & $o(h^{3/2-\delta})$ \\
digital   & $O(h^{1/2})$ & $O(h^{1/2}\log h)$ & $O(h^{3/2})$ & $o(h^{3/2-\delta})$ \\ \hline
\end{tabular}
\end{center}
\caption{Observed and theoretical convergence rates for the multilevel 
correction variance for scalar SDEs, using the Euler-Maruyama and Milstein 
discretisations. $\delta$ is any strictly positive constant.}
\label{tab:EM-Milstein}
\end{table}

\subsection{Milstein discretisation}

For Lipschitz payoffs, the variance $V_\ell$ for the natural 
multilevel estimator converges at twice the order of the strong 
convergence of the numerical approximation of the SDE.  This
immediately suggests that it would be better to replace the 
Euler-Maruyama discretisation by the Milstein discretisation
\cite{giles08b}
since it gives first order strong convergence under certain 
conditions (see Theorem 10.3.5 in \cite{kp92}).

This immediately gives an improved variance for European and Asian 
options, as shown in Table \ref{tab:EM-Milstein}, but to get the 
improved variance for lookback, barrier and digital options requires
the construction of estimators which are slightly different on the 
coarse and fine path simulations, but which respect the condition that
$\mathbb{E}[P^f_\ell] = \mathbb{E}[P^c_\ell]$.  

The construction for the digital option will be discussed next, but 
for the lookback and barrier options, the key is the definition of a 
Brownian Bridge interpolant based on the approximation that the drift 
and volatility do not vary within the timestep. For each coarse 
timestep, the mid-point of the interpolant can be sampled using 
knowledge of the fine path Brownian increments, and then classical 
results can be used for the distribution of the minimum or maximum within 
each fine timestep for both the fine and coarse path approximations 
\cite{glasserman04}.  The full details are given in \cite{giles08b},
and Table \ref{tab:EM-Milstein} summarises the convergence behaviour
observed numerically, and the supporting numerical analysis by
Giles, Debrabant \& R\"{o}{\ss}ler \cite{gdr13}.

The outcome is that for the case in which the number of timesteps 
doubles at each level, so $h_\ell = 2^{-\ell} h_0$, then 
$\gamma\!=\!1$ and either $\beta\!=\!2$ (European, Asian and lookback) 
or $\beta\!=\!1.5$ (barrier and digital). Hence, we are in the regime 
where $\beta\!>\!\gamma$ and the overall complexity is $O(\varepsilon^{-2})$.
Furthermore, the dominant computational cost is on the coarsest 
levels of simulation.  

Since the coarsest levels are low-dimensional, they are well
suited to the use of quasi-Monte Carlo methods which are 
particularly effective in lower dimensions because of the 
existence of $O((\log N)^d/N)$ error bounds, where
$d$ is the dimension and $N$ is the number of QMC points.  
The bounds are for the numerical integration of certain function
classes on the unit hypercube, and are a consequence of the 
Koksma-Hlawka inequality together with bounds on the 
star-discrepancy of certain sequences of QMC points.

This has been investigated by Giles \& Waterhouse \cite{gw09}
using a rank-1 lattice rule to generate the quasi-random 
numbers, randomisation with 32 independent offsets to obtain 
confidence intervals, and a standard Brownian Bridge construction 
of the increments of the driving Brownian process.
The numerical results show that MLMC on its own was better 
than QMC on its own, but the combination of the two was even better.
The QMC treatment greatly reduced the variance per sample for the 
coarsest levels, resulting in significantly reduced costs overall.  
In the simplest case of a Lipschitz European payoff, 
the computational complexity was reduced from $O(\varepsilon^{-2})$ to 
approximately $O(\varepsilon^{-1.5})$.  


\subsubsection{Digital options}

As discussed earlier, discontinuous payoffs pose a challenge to the 
multilevel Monte Carlo approach, because small differences in the 
coarse and fine path simulations can lead to an $O(1)$ difference in 
the payoff function.  This leads to a slower decay in the variance 
$V_\ell$, and because the fourth moment is also much larger it leads 
to more samples being required to obtain an accurate estimate for 
$V_\ell$, which is needed to determine the optimal number of samples 
$N_\ell$.  

This is a generic problem.  Although we will discuss it 
here in the specific context of a Brownian SDE and an option which 
is a discontinuous function of the underlying at the final time, 
the methods which are discussed are equally applicable in a range 
of other cases.  Indeed, some of these techniques have been first 
explored in the context of pathwise sensitivity analysis \cite{bg12}
or jump-diffusion modelling \cite{xg12}.

\vspace{0.1in}
{\bf Conditional expectation}

The conditional expectation approach builds on a well-established
technique for payoff smoothing which is used for pathwise sensitivity 
analysis (see, for example, pp.~399-400 in \cite{glasserman04}).

We start by considering the fine path simulation, and make a slight 
change by using the Euler-Maruyama discretisation for the final timestep,
instead of the Milstein discretisation. Conditional on the numerical
approximation of the value $S_{T\!-h}$ one timestep before the end
(which in turn depends on all of the Brownian increments up to that time) 
the numerical approximation for the final value $S_T$ now has a Gaussian 
distribution, and for a simple digital option the conditional expectation
is known analytically.

The same treatment is used for the coarse path, except that in the 
final timestep, we re-use the known value of the Brownian increment 
for the second last fine timestep, which corresponds to the first half
of the final coarse timestep. This results in the conditional distribution 
for the coarse path underlying at maturity matching that of the fine 
path to within $O(h)$, for both the mean and the standard deviation
\cite{gdr13}.  Consequently, the difference in payoff between the 
coarse and fine paths near the payoff discontinuity is $O(h^{1/2})$,
and so the variance is approximately $O(h^{3/2})$.

\vspace{0.1in}
{\bf Splitting}

The conditional expectation technique works well in 1D where there 
is a known analytic value for the conditional expectation, but in 
multiple dimensions it may not be known.  In this case, one can use 
the technique of ``splitting'' \cite{ag07}. Here the conditional 
expectation is replaced by a numerical estimate, averaging over a 
number of sub-samples. i.e.~for each set of Brownian increments up
to one fine timestep before the end, one uses a number of samples 
of the final Brownian increment to produce an average payoff.  If
the number of sub-samples is chosen appropriately, the variance 
is the same, to leading order, without any increase in the 
computational cost, again to leading order. Because of its simplicity 
and generality, this is now my preferred approach.  Furthermore, 
one can revert to using the Milstein approximation for the final 
timestep.

\vspace{0.1in}
{\bf Change of measure}

The change of measure approach is another approximation to the 
conditional expectation.  The fine and coarse path conditional 
distributions at maturity are two very similar Gaussian distributions.
Instead of following the splitting approach of taking corresponding 
samples from these two distributions, we can instead take a sample
from a third Gaussian distribution (with a mean and variance perhaps
equal to the average of the other two).  This leads to the introduction
of a Radon-Nikodym derivative for each path, and the difference in the
payoffs from the two paths is then due to the difference in their 
Radon-Nikodym derivatives.

In the specific context of digital options, this is a more complicated
method to implement, and the resulting variance is no better.  However,
in other contexts a similar approach can be very effective.

\subsubsection{Multi-dimensional SDEs}

The discussion so far has been for scalar SDEs, but the computational 
benefits of Monte Carlo methods arise in higher dimensions.  For 
multi-dimensional SDEs satisfying the usual commutativity condition 
(see, for example, p.353 in \cite{glasserman04}) the Milstein
discretisation requires only Brownian increments for its implementation,
and most of the analysis above carries over very naturally.  

The only difficulties are in lookback and barrier options where the 
classical results for the distribution of the minimum or maximum of a 
one-dimensional Brownian motion, do not extend to the joint distribution 
of the minima or maxima of two correlated Brownian motions. An alternative
approach may be to sub-sample from the Brownian Bridge interpolant for
those timesteps which are most likely to give the global minimum or 
maximum.  This may need to be combined with splitting for the barrier 
option to avoid the $O(1)$ difference in payoffs.  An alternative might 
be to use adaptive time-stepping \cite{hsst12}.

For multi-dimensional SDEs which do not satisfy the commutativity 
condition the Milstein discretisation requires the simulation of 
L{\'e}vy areas.  This is unavoidable to achieve first order strong 
convergence; the classical result of Clark \& Cameron says that 
$O(h^{1/2})$ strong convergence is the best that can be achieved 
in general using just Brownian increments \cite{cc80}.

However, Giles \& Lukasz have developed an antithetic treatment
which achieves a very low variance despite the $O(h^{1/2})$ 
strong convergence \cite{gs12}.  The estimator which is used is
\[
Y_\ell = N_\ell^{-1}\ \sum_i\ 
{\textstyle \frac{1}{2}} \left( P_\ell(\omega_i) \!+\! P_\ell(\omega^a_i) \right)
- P_{\ell-1}(\omega_i).
\]
Here $\omega_i$ represents the driving Brownian path, and $\omega^a_i$
is an antithetic counterpart defined by a time-reversal of the Brownian 
path within each coarse timestep.  This results in the Brownian increments
for the antithetic fine path being swapped relative to the original path.
Lengthy analysis proves that the average of the fine and antithetic paths 
is within $O(h)$ of the coarse path, and hence the multilevel variance 
is $O(h^2)$ for smooth payoffs, and $O(h^{3/2})$ for the standard European 
call option.

This treatment has been extended to handle lookback and barrier options 
\cite{gs13}.  This combines sub-sampling of the Brownian path to 
approximate the L{\'e}vy areas with sufficient accuracy to achieve 
$O(h^{3/4})$ strong convergence, with an antithetic treatment at the 
finest level of resolution to ensure that the average of the fine paths 
is within $O(h)$ of the coarse path.

\subsection{L{\'e}vy processes}

\subsubsection{Jump-diffusion processes}

With finite activity jump-diffusion processes, such as in the Merton
model \cite{merton76}, it is natural to simulate each individual 
jump using a jump-adapted discretisation \cite{pb10}.

If the jump rate is constant, then the jumps on the coarse and 
fine paths will occur at the same time, and the extension of the 
multilevel method is straightforward \cite{xg12}.

If the jump rate is path-dependent then the situation is trickier.
If there is a known upper bound to the jump rate, then one 
can use Glasserman \& Merener's ``thinning'' approach \cite{gm04}
in which a set of candidate jump times is simulated based on 
the constant upper bound, and then a subset of these are selected 
to be real jumps.  The problem with the multilevel extension of
this is that some candidate jumps will be selected for the coarse 
path but not for the fine path, or vice versa, leading to an $O(1)$
difference in the paths and hence the payoffs.  Xia overcomes this 
by using a change of measure to select the jump times consistently
for both paths, with a Radon-Nikodym derivative being introduced 
in the process \cite{xg12}.

\subsubsection{More general processes}

With infinite activity L{\'e}vy processes it is impossible to 
simulate each jump.
One approach is to simulate the large jumps and either neglect 
the small jumps or approximate their effect by adding a Brownian 
diffusion term \cite{dereich11,dh11,marxen10}.  
Following this approach, the cutoff $\delta_\ell$ 
for the jumps which are simulated varies with level, and
$\delta_\ell\rightarrow 0$ as $\ell\rightarrow \infty $
to ensure that the bias converges to zero. In the multilevel
treatment, when simulating $P_\ell - P_{\ell-1}$ the jumps 
fall into three categories.  The ones which are larger than 
$\delta_{\ell-1}$ get simulated in both the fine and coarse paths.
The ones which are smaller than $\delta_\ell$ are either 
neglected for both paths, or approximated by the same Brownian 
increment.  The difficulty is in the intermediate range 
$[\delta_\ell, \delta_{\ell-1}]$ in which the jumps are 
simulated for the fine path, but neglected or approximated 
for the coarse path.  This is what leads to the difference 
in path simulations, and hence to a non-zero value for
$P_\ell - P_{\ell-1}$.

Alternatively, for many SDEs driven by a L{\'e}vy process it is
possible to directly simulate the increments of the L{\'e}vy 
process over a set of uniform timesteps \cite{ct04,schoutens03}, 
in exactly the same way as one simulates Brownian increments.
For other L{\'e}vy processes, it may be possible in the future
to simulate the increments by constructing approximations to the 
inverse of the cumulative distribution function. 
Where this is possible, it may be the best approach to achieve
a close coupling between the coarse and fine path simulations,
and hence a low variance $V_\ell$,  
since the increments of the driving L{\'e}vy process for the coarse 
path can be obtained trivially by summing the increments for the
fine path.


\section{SPDEs}

After developing the MLMC method for SDE simulations, it was 
immediately clear that it was equally applicable to SPDEs, and 
indeed the computational savings would be greater because the
cost of a single sample increases more rapidly with grid
resolution for SPDEs with higher space-time dimension.  

In 2006, the author discussed this with Thomas Hou in the 
specific context of elliptic SPDEs with random coefficients, 
and Hou's postdoc then performed the first unpublished MLMC 
computations for SPDEs.
The first published work was by a student of Klaus Ritter
in her Diploma thesis \cite{graubner08}; 
her application was to parabolic SPDEs.  
Since this early work, there has been a variety of papers on 
elliptic \cite{bsz11,cst13,cgst11,tsgu13}, 
parabolic \cite{bl12,gr12}
and hyperbolic \cite{mss12} SPDEs.

In almost all of this work, the construction of the multilevel 
estimator is quite natural, using a geometric sequence of grids
and the usual estimators for $P_\ell \!-\!P_{\ell-1}$.
It is the numerical analysis of the 
variance of the multilevel estimator which is often very challenging.

\subsection{Elliptic SPDE}

The largest amount of research on multilevel for SPDEs has been for
elliptic PDEs with random coefficients.  The PDE typically has the form
\begin{equation*}
- \nabla \cdot \left(k(\mathbf{x},\omega) \nabla p(\mathbf{x}, \omega) \right) 
= 0, \qquad \mathbf{x} \in D.
\end{equation*}
with Dirichlet or Neumann boundary conditions on the boundary $\partial D$.
For sub-surface flow problems, such as the modelling of groundwater flow in 
nuclear waste repositories, the diffusivity (or permeability) $k$ is often 
modelled as a lognormal random field, i.e.~$\log k$ is a Gaussian field with 
a uniform mean (which we will take to be zero for simplicity) and a covariance 
function of the general form 
$R(\mathbf{x}, \mathbf{y}) = r(\mathbf{x}\!-\!\mathbf{y})$.
Samples of $\log k$ are provided by a Karhunen-Lo\`eve expansion:
\[
\log k(\mathbf{x},\omega) = 
\sum_{n=0}^{\infty}\sqrt{\theta_n}\ \xi_n(\omega)\ f_n(\mathbf{x}),
\]
where $\theta_n$ are the eigenvalues of $R(\mathbf{x}, \mathbf{y})$ in
decreasing order, $f_n$ are the corresponding eigenfunctions, 
and $\xi_n$ are independent unit Normal random variables.  However, 
it is more efficient to generate them using a circulant embedding 
technique which enables the use of FFTs \cite{dn97}.

The multilevel treatment is straightforward.  The spatial grid resolution
is doubled on each level.  Using the Karhunen-Lo\`eve generation, the 
expansion is truncated after $K_\ell$ terms, with $K_\ell$ increasing 
with level \cite{tsgu13}; in unpublished work, a similar approach has 
also been used with the circulant embedding generation.  

In both cases, 
$\log k$ is generated using a row-vector of independent unit Normal 
random variables $\xi$.  The variables for the fine level can be 
partitioned into those for the coarse level $\xi_{\ell\!-\!1}$, 
plus some additional variables $z_\ell$, 
giving $\xi_\ell = (\xi_{\ell-1}, z_\ell)$.
It is possible to develop an antithetic treatment similar to that 
used for SDEs by defining $\xi^a_\ell = (\xi_{\ell-1}, -z_\ell)$.  
This gives a second $\log k_\ell^a$ field on the fine grid, and 
then the multilevel estimator can be based on the average of the 
two outputs obtained on the fine grid, minus the output obtained 
on the coarse grid using $\log k_{\ell-1}$.  Unfortunately, numerical
experiments indicate it gives little benefit; it is mentioned here
as another illustration of an antithetic estimator, and as a 
warning that it does not always yields significant benefits.

The numerical analysis of the multilevel approach for these elliptic SPDE
applications is challenging because the diffusivity is unbounded, but
Charrier, Scheichl \& Teckentrup \cite{cst13} have successfully analysed 
it for certain output functionals, and Teckentrup {\it et al} have 
further developed the analysis for other output functionals and more
general log-normal diffusivity fields \cite{tsgu13}.

\subsection{Parabolic SPDE}

Giles \& Reisinger \cite{gr12} consider an unusual SPDE 
from credit default modelling,
\[
\D p = -\mu\, \frac{\partial p}{\partial x}\ \D t 
+ \frac{1}{2} \, \frac{\partial^2  p }{\partial x^2}\ \D t
- \sqrt{\rho}\ \frac{\partial p}{\partial x}\ \D M_t,
~~~ x>0
\]
subject to boundary condition $p(0,t) \!=\!0$.
Here $p(x,t)$ represents the probability density function for firms 
being a distance $x$ from default at time $t$. The diffusive term 
is due to idiosyncratic factors affecting individual firms, while 
the stochastic term due to the scalar Brownian motion $M_t$ 
corresponds to the systemic movement due to random market effects 
affecting all firms. 
The payoff corresponds to different tranches of a credit derivative
which depends on the integral
$\int_0^\infty p(x,t) \ \D x$ at a set of discrete times.

A Milstein time discretisation with timestep $k$, and a central 
space discretisation of the spatial derivatives with uniform 
spacing $h$ gives the numerical approximation
\[
p_j^{n+1} = p_j^n\ -\ \frac{\mu\, k + \sqrt{\rho\, k}\, Z_n}{2h} 
\left(p_{j+1}^n - p_{j-1}^n\right) 
\nonumber + \frac{(1\!-\!\rho)\, k + \rho \, k\, Z_n^2}{2h^2} 
 \left(p_{j+1}^n - 2 p_j^n + p_{j-1}^n\right)
\]
where $p_j^n \approx p(j\, h, n\, k)$, 
and the $Z_n$ are standard Normal random variables so that
$\sqrt{h}\ Z_n$ corresponds to an increment of the driving scalar
Brownian motion.

The multilevel implementation is very straightforward, with 
$k_\ell = k_{\ell-1}/2$ and $h_\ell = h_{\ell-1}/4$ due to 
numerical stability considerations which are analysed in the paper.
As with SDEs, the coupling between the coarse and fine samples 
comes from summing the fine path Brownian increments in pairs to
give the increments for the coarse path.  
The computational cost increases by factor 8 on each level, and 
numerical experiments indicate that the variance decreases by 
factor 8, so the overall computational complexity to achieve an 
$O(\varepsilon)$ RMS error is again $O(\varepsilon^{-2} (\log \varepsilon)^2)$.

\section{Continuous-time Markov Chain simulation}

Anderson \& Higham have recently developed a very interesting new
application of multilevel to continuous-time Markov Chain simulation
\cite{ah12}.  Although they present their work in the context of
stochastic chemical reactions, when species concentrations are 
extremely low and so stochastic effects become significant, they
point out that the method has wide applicability in other areas.

In the simplest case of a single chemical reaction, the 
``tau-leaping'' method (which is essentially the Euler-Maruyama 
method, approximating the reaction rate as being constant throughout
the timestep) gives the discrete equation 
\[
{\bf x}_{n+1} = {\bf x}_{n} + P( h\ \lambda({\bf x}_{n})),
\]
where $h$ is the timestep, $\lambda({\bf x}_{n})$ is the reaction 
rate (or propensity function), and $P(t)$
represents a unit-rate Poisson random variable over time interval $t$.

If this equation defines the fine path in the multilevel simulation,
then the coarse path, with double the timestep, is given by
\[
{\bf x}^c_{n+2} = {\bf x}^c_{n} + P( 2h \ \lambda({\bf x}^c_{n}))
\]
for even timesteps $n$.  The question then is how to couple the
coarse and fine path simulations.

The key observation by Anderson \& Higham \cite{ah12} 
is that for any $t_1, t_2 > 0$, 
the sum of two independent Poisson variates $P(t_1)$, $P(t_2)$
is equivalent in distribution to $P(t_1 \!+\! t_2)$.
Based on this, the first step is to express the coarse path 
Poisson variate as the sum of two Poisson variates, 
$P( h \, \lambda({\bf x}^c_{n}))$ corresponding to the first 
and second fine path timesteps.
For the first of the two fine timesteps, the coarse and fine path 
Poisson variates are coupled by defining two Poisson variates
based on the minimum of the two reactions rates, and the absolute 
difference,
\[
P_1 = P\left(\rule{0in}{0.14in} h 
\min(\lambda({\bf x}_{n}),\lambda({\bf x}^c_{n})) \right),
~~~
P_2 = P\left(\rule{0in}{0.14in} h 
\left| \lambda({\bf x}_{n})-\lambda({\bf x}^c_{n}) \right| \right),
\]
and then using $P_1$ as the Poisson variate for the path with the 
smaller rate, and $P_1\!+\!P_2$ for the path with the larger rate.
This elegant approach naturally gives a small difference in the 
Poisson variates when the difference in rates is small, and leads 
to a very effective multilevel algorithm.

In their paper \cite{ah12}, Anderson \& Higham treat more general
systems with multiple reactions, and include an additional coupling 
at the finest level to an SSA (Stochastic Simulation Algorithm) 
computation, so that their overall multilevel estimator is unbiased,
unlike the estimators discussed earlier for SDEs.
Finally, they give a complete numerical analysis of
the variance of their multilevel algorithm.

Because stochastic chemical simulations typically involve 1000's 
of reactions, the multilevel method is particularly effective in
this context, providing computational savings in excess of a 
factor of 100 \cite{ah12}.

\section{Wasserstein metric}

In the multilevel treatment of SDEs, the Brownian or L{\'e}vy  
increments for the coarse path are obtained by summing the
increments for the fine path.  Similarly, in the Markov Chain
treatment, the Poisson variate for the coarse timestep is defined 
as the sum of two Poisson variates for fine timesteps.  

This sub-division of coarse path random variable into the sum 
of two fine path random variables should work in many settings.
The harder step in more general applications is likely to be the 
second step in the Markov Chain treatment, tightly coupling
the increments used for the fine and coarse paths over the same 
fine timestep.

The general statement of this problem is the following: given
two very similar scalar probability distributions, we want to 
obtain samples $Z_f, Z_c$ from each in a way which minimises 
$\mathbb{E}[\,|Z_f\!-\!Z_c|^p]$.  This corresponds precisely to the 
Wasserstein metric which defines the ``distance'' between 
two probability distributions as
\[
\left( 
\inf_\gamma \int \left\|Z_f \!-\! Z_c\right\|^p \, d\gamma(Z_f, Z_c)
\right)^{1/p},
\]
where the minimum is over all joint distributions with the 
correct marginals.  In 1D, the Wasserstein metric is equal to
\[
\left( 
\int_0^1 \left| \Phi_f^{-1}(u) -\Phi_c^{-1}(u) \right|^p \, \D u
\right)^{1/p},\]
where $\Phi_f$ and $\Phi_c$ are the cumulative probability 
distributions for $Z_f$ and $Z_c$ \cite{bf81}, 
and this minimum is achieved by choosing
$Z_f = \Phi_f^{-1}(U)$, $Z_c = \Phi_c^{-1}(U)$,
for the same uniform $[0,1]$ random variable $U$.
This suggests this may be a good general technique for future 
multilevel applications, provided one is able to invert the 
relevant cumulative distributions, possibly through generating 
appropriate spline approximations.

\section{Other uses of multilevel}

\subsection{Nested simulation}

The pricing of American options is one of the big challenges for
Monte Carlo methods in computational finance, and
Belomestny \& Schoenmakers have recently written a very interesting
paper on the use of multilevel Monte Carlo for this purpose
\cite{bs13}. Their method is based on Anderson \& Broadie's
dual simulation method \cite{ab04} in which a key component at each 
timestep in the simulation is to estimate a conditional expectation 
using a number of sub-paths.

In their multilevel treatment, Belomestny \& Schoenmakers use the
same uniform timestep on all levels of the simulation.  The quantity 
which changes between different levels of simulation is the number of
sub-samples used to estimate the conditional expectation.  To couple 
the coarse and fine levels, the fine level uses $N_\ell$ sub-samples, 
and the coarse level uses $N_{\ell-1} = N_\ell/2$ of them. 

Related unpublished research by N.~Chen
for a similar multilevel treatment of nested simulation found 
that the multilevel correction variance is reduced if the payoff 
on the coarse level is replaced by an average of the payoffs obtained 
using the first $N_\ell/2$ and the second $N_\ell/2$ samples.
This is similar in some ways to the antithetic approach described 
earlier.

In future research, Belomestny \& Schoenmakers intend to also 
change the number of timesteps on each level, to increase the 
overall computational benefits of the multilevel approach.

\subsection{Truncated series expansions}

Building on earlier work by Broadie \& Kaya \cite{bk06},
Glasserman \& Kim have recently developed an efficient 
method \cite{gk11} of exactly simulating the Heston 
stochastic volatility model \cite{heston93}.
The key to their algorithm is a method of representing
the integrated volatility over a time interval $[0,T]$,
conditional on the initial and final values, $v_0$ and $v_T$ as 
\[
\left(\left.\int_0^T V_s \, ds\ \right|\ V_0=v_0, V_T= v_T\right)
\ \stackrel{d}{=}\ 
\sum_{n=1}^\infty x_n + 
\sum_{n=1}^\infty y_n + 
\sum_{n=1}^\infty z_n
\]
where $x_n, y_n, z_n$ are independent random variables.

In practice, they truncate the series expansions at a level which 
ensures the desired accuracy, but a more severe truncation would
lead to a tradeoff between accuracy and computational cost.
This makes the algorithm a candidate for a multilevel treatment
in which the level $\ell$ computation performs the truncation at 
$N_\ell$, so the level $\ell$ computation would use
\[
\sum_{n=1}^{N_\ell} x_n + 
\sum_{n=1}^{N_\ell} y_n + 
\sum_{n=1}^{N_\ell} z_n
\]
while the level $\ell \!-\! 1$ computation would truncate the 
summations at $N_{\ell- 1}$, but would use the same random variables 
$x_n, y_n, z_n$ for $1\leq n \leq N_{\ell- 1}$.

This kind of multilevel treatment has not been tested experimentally, 
but it seems that it might yield some computational savings even though
Glasserman \& Kim typically only need to retain 10 terms in their 
summations through the use of a carefully constructed estimator for the 
truncated remainder.  The savings may be larger in other circumstances 
which require more terms to be retained for the desired accuracy.

\subsection{Mixed precision arithmetic}

The final example of the use of multilevel is unusual, because it 
concerns the computer implementation of Monte Carlo algorithms.
In the latest CPUs from Intel and AMD, each core has a vector
unit which can perform 8 single precision or 4 double precision
operations with one instruction.  Also, double precision data 
takes twice as much time to transfer as single precision data.
Hence, single precision computations can be twice as fast as 
double precision on CPUs, and the difference can be even greater 
on GPUs.  This raises the question of whether single precision 
arithmetic is sufficient for Monte Carlo simulation. 

My view is that it usually is since the finite precision rounding 
errors are smaller than the other sources of error: statistical 
error due to Monte Carlo sampling; bias due to SDE discretisation;
model uncertainty.  However, there can be significant errors when 
averaging unless one uses binary tree summation \cite{higham93} 
to perform the summation, and in addition computing sensitivities 
by perturbing input parameters (so-called ``bumping'') can greatly 
amplify the rounding errors.

The best solution is perhaps to use double precision for the final 
averaging, and pathwise sensitivity analysis or the likelihood ratio 
method for computing sensitivities, but if there remains a need 
for the path simulation to be performed in double precision then 
one could use the two-level MLMC approach in which level 0 corresponds 
to single precision and level 1 corresponds to double precision,
with the same random numbers being used for both.

\subsection{Multiple outputs}

In all of the discussion so far, we have been concerned with a single
expectation arising from a stochastic simulation.  However, there are
often times when one wishes to estimate the expected value of multiple 
outputs.

Extending the analysis in section \ref{sec:unbiased_MLMC}, when using 
multilevel 
to estimate $M$ different expectations, using $N_l$ samples on each level, 
the goal is to achieve an acceptably small variance for each output
\[
\sum_{\ell=0}^L N_\ell^{-1}\, V_{\ell,m}\ \leq\ \varepsilon_m^2, 
~~~~~ m=1, \ldots, M,
\]
with the desired accuracy $\varepsilon_m$ being allowed to vary from one output
to another, and to do so with the minimum computational cost which is 
given as usual as
\[
\sum_{\ell=0}^L N_\ell \, C_\ell,
\]
assuming that the cost of computing the output functions is negligible
compared to the cost of obtaining the stochastic sample (e.g.~through an 
SDE path simulation). 

This leads naturally to a constrained optimisation problem with a 
separate Lagrange multiplier for each output.  However, a much 
simpler idea, due to Tigran Nagapetyan, 
which in practice is almost always equivalent, is to define
\[
V_\ell = \max_m \frac{V_{\ell,m}}{\varepsilon_m^2}
\]
and make the variance constraint
$\displaystyle
\ \sum_{\ell=0}^L N_\ell^{-1}\, V_\ell\ \leq\ 1.
$

This is sufficient to ensure that all of the individual constraints are 
satisfied, and we can then use the standard approach with a single Lagrange
multiplier.
%
This multi-output approach is currently being investigated by Nagapetyan, 
Ritter and the author for the approximation of cumulative distribution 
functions and probability density functions arising from stochastic 
simulations.


\section{Conclusions}

In the past 6 years, considerable progress has been achieved with 
the multilevel Monte Carlo method for a wide range of applications.
This review has attempted to emphasise the conceptual simplicity 
of the multilevel approach; in essence it is simply a recursive 
control variate strategy, using cheap approximations to some random 
output quantity as a control variate for more accurate but more costly 
approximations.  

In practice, the challenge is to develop a tight coupling between 
successive approximation levels, to minimise the variance of the 
difference in the output obtained from each level. In the context 
of SDE and SPDE simulations, strong convergence properties are often 
relied on to obtain a small variance between coarse and fine simulations.
In the specific context of a digital option associated with a Brownian 
SDE, three treatments were described to effectively smooth the output:
a analytic conditional expectation, a ``splitting'' approximation, and
a change of measure.  Similar treatments have been found to be helpful 
in other contexts.

Overall, multilevel methods are being used for an increasingly wide 
range of applications.  The biggest savings are in situations in 
which the coarsest approximation is very much cheaper than the finest.
So far, this includes multi-dimensional SPDEs, and chemical stochastic
simulations with 1000's of timesteps.  In SDE simulations which perhaps 
only require 32 timesteps for the desired level of accuracy, the
potential savings are naturally quite limited.

\vspace{0.05in}

Although this is primarily a survey article, a few new ideas have 
been introduced:
\begin{itemize}
\item 
equation (\ref{eq:total_cost}) giving the total computational 
cost required for a general unbiased multilevel estimator is new, 
as is the discussion which follows it, although the underlying 
analysis is not;

\item
based on the 1D Wasserstein metric, it seems that inverting the 
relevant cumulative distributions may be a good way to couple 
fine and coarse level simulations in multilevel implementations;

\item
the multilevel approach could be used in applications which involve 
the truncation of series expansions;

\item
a two-level method combining single and double precision 
computations might provide useful savings, due to the 
lower cost of single precision arithmetic;

\item
a multilevel approach for situations with multiple expectations to be estimated.
\end{itemize}

Looking to the future, exciting areas for further research include:
\begin{itemize}
\item
more use of multilevel for nested simulations;

\item
further investigation of multilevel quasi-Monte Carlo methods;

\item
continued research on numerical analysis, especially for SPDEs;

\item
development of multilevel estimators for new applications.
\end{itemize}


For further information on multilevel Monte Carlo methods, see the webpage\\
{\tt http://people.maths.ox.ac.uk/gilesm/mlmc\_community.html}\\
which lists the research groups working in the area, and their 
main publications.


\begin{thebibliography}{10}

\bibitem{ab04}
L.~Andersen and M.~Broadie.
\newblock A primal-dual simulation algorithm for pricing multi-dimensional
  {A}merican options.
\newblock {\em Management Science}, 50(9):1222--1234, 2004.

\bibitem{ah12}
D.~Anderson and D.J. Higham.
\newblock Multi-level {M}onte {C}arlo for continuous time {M}arkov chains with
  applications in biochemical kinetics.
\newblock {\em SIAM Multiscale Modeling and Simulation}, 10(1):146--179, 2012.

\bibitem{ag07}
A.~Asmussen and P.~Glynn.
\newblock {\em Stochastic Simulation}.
\newblock Springer, New York, 2007.

\bibitem{avikainen09}
R.~Avikainen.
\newblock On irregular functionals of {SDEs} and the {E}uler scheme.
\newblock {\em Finance and Stochastics}, 13(3):381--401, 2009.

\bibitem{bl12}
A.~Barth and A.~Lang.
\newblock Multilevel {M}onte {C}arlo method with applications to stochastic
  partial differential equations.
\newblock {\em Int.~Journal of Computer Mathematics}, 89(18):2479--2498, 2012.

\bibitem{bsz11}
A.~Barth, C.~Schwab, and N.~Zollinger.
\newblock Multi-level {M}onte {C}arlo finite element method for elliptic {PDE}s
  with stochastic coefficients.
\newblock {\em Numerische Mathematik}, 119(1):123--161, 2011.

\bibitem{bs13}
D.~Belomestny and J.~Schoenmakers.
\newblock Multilevel dual approach for pricing {A}merican style derivatives.
\newblock {\em Finance and Stochastics}, 2013.

\bibitem{bf81}
P.J. Bickel and D.A. Freedman.
\newblock Some asymptotic theory for the bootstrap.
\newblock {\em Annals of Statistics}, 9:1196--1217, 1981.

\bibitem{bgr94}
A.~Brandt, M.~Galun, and D.~Ron.
\newblock Optimal multigrid algorithms for calculating thermodynamic limits.
\newblock {\em Journal of Statistical Physics}, 74(1-2):313--348, 1994.

\bibitem{bi03}
A.~Brandt and V.~Ilyin.
\newblock Multilevel {M}onte {C}arlo methods for studying large scale phenomena
  in fluids.
\newblock {\em Journal of Molecular Liquids}, 105(2-3):245--248, 2003.

\bibitem{bk06}
M.~Broadie and O.~Kaya.
\newblock Exact simulation of stochastic volatility and other affine jump
  diffusion processes.
\newblock {\em Operations Research}, 54(2):217--231, 2006.

\bibitem{bg12}
S.~Burgos and M.B. Giles.
\newblock Computing {G}reeks using multilevel path simulation.
\newblock In L.~Plaskota and H.~Wo{\'z}niakowski, editors, {\em Monte Carlo and
  Quasi-Monte Carlo Methods 2010}, pages 281--296. Springer-Verlag, 2012.

\bibitem{cst13}
J.~Charrier, R.~Scheichl, and A.~Teckentrup.
\newblock Finite element error analysis of elliptic {PDE}s with random
  coefficients and its application to multilevel {M}onte {C}arlo methods.
\newblock {\em SIAM Journal on Numerical Analysis}, 51(1):322--352, 2013.

\bibitem{cc80}
J.M.C. Clark and R.J. Cameron.
\newblock The maximum rate of convergence of discrete approximations for
  stochastic differential equations.
\newblock In B.~Grigelionis, editor, {\em Stochastic Differential Equations},
  volume~25 of {\em Lecture Notes in Control and Information Sciences}.
  Springer-Verlag, 1980.

\bibitem{cgst11}
K.A. Cliffe, M.B. Giles, R.~Scheichl, and A.~Teckentrup.
\newblock Multilevel {M}onte {C}arlo methods and applications to elliptic
  {PDE}s with random coefficients.
\newblock {\em Computing and Visualization in Science}, 14(1):3--15, 2011.

\bibitem{ct04}
R.~Cont and P.~Tankov.
\newblock {\em Financial modelling with jump processes}.
\newblock CRC Press, 2004.

\bibitem{dereich11}
S.~Dereich.
\newblock Multilevel {M}onte {C}arlo algorithms for {L}{\'e}vy-driven {SDE}s
  with {G}aussian correction.
\newblock {\em Annals of Applied Probability}, 21(1):283--311, 2011.

\bibitem{dh11}
S.~Dereich and F.~Heidenreich.
\newblock A multilevel {M}onte {C}arlo algorithm for {L}{\'e}vy-driven
  stochastic differential equations.
\newblock {\em Stochastic Processes and their Applications}, 121(7):1565--1587,
  2011.

\bibitem{dn97}
C.R. Dietrich and G.H. Newsam.
\newblock Fast and exact simulation of stationary {G}aussian processes through
  circulant embedding of the covariance matrix.
\newblock {\em SIAM Journal on Scientific Computing}, 18(4):1088–1107, 1997.

\bibitem{giles08b}
M.B. Giles.
\newblock Improved multilevel {M}onte {C}arlo convergence using the {M}ilstein
  scheme.
\newblock In A.~Keller, S.~Heinrich, and H.~Niederreiter, editors, {\em Monte
  Carlo and Quasi-Monte Carlo Methods 2006}, pages 343--358. Springer-Verlag,
  2008.

\bibitem{giles08}
M.B. Giles.
\newblock Multilevel {M}onte {C}arlo path simulation.
\newblock {\em Operations Research}, 56(3):607--617, 2008.

\bibitem{giles09}
M.B. Giles.
\newblock Multilevel {M}onte {C}arlo for basket options.
\newblock In M.D. Rossetti, R.R. Hill, B.~Johansson, A.~Dunkin, and R.G.
  Ingalls, editors, {\em Proceedings of the 2009 Winter Simulation Conference},
  pages 1283--1290. IEEE, 2009.

\bibitem{gdr13}
M.B. Giles, K.~Debrabant, and A.~R{\"o}{\ss}ler.
\newblock Numerical analysis of multilevel {M}onte {C}arlo path simulation
  using the {M}ilstein discretisation.
\newblock {\em ArXiv preprint: 1302.4676}, 2013.

\bibitem{ghm09}
M.B. Giles, D.J. Higham, and X.~Mao.
\newblock Analysing multilevel {M}onte {C}arlo for options with non-globally
  {L}ipschitz payoff.
\newblock {\em Finance and Stochastics}, 13(3):403--413, 2009.

\bibitem{gr12}
M.B. Giles and C.~Reisinger.
\newblock Stochastic finite differences and multilevel {M}onte {C}arlo for a
  class of {SPDEs} in finance.
\newblock {\em SIAM Journal of Financial Mathematics}, 3(1):572--592, 2012.

\bibitem{gs12}
M.B. Giles and L.~Szpruch.
\newblock Antithetic multilevel {M}onte {C}arlo estimation for
  multi-dimensional {SDE}s without {L\'e}vy area simulation.
\newblock {\em ArXiv preprint: 1202.6283}, 2012.

\bibitem{gs13}
M.B. Giles and L.~Szpruch.
\newblock Antithetic multilevel {M}onte {C}arlo estimation for multidimensional
  {SDE}s.
\newblock In {\em Monte Carlo and Quasi-Monte Carlo Methods 2012 (submitted)}.
  Springer-Verlag, 2013.

\bibitem{gw09}
M.B. Giles and B.J. Waterhouse.
\newblock Multilevel quasi-{M}onte {C}arlo path simulation.
\newblock In {\em Advanced Financial Modelling}, Radon Series on Computational
  and Applied Mathematics, pages 165--181. De Gruyter, 2009.

\bibitem{glasserman04}
P.~Glasserman.
\newblock {\em {M}onte {C}arlo Methods in Financial Engineering}.
\newblock Springer, New York, 2004.

\bibitem{gk11}
P.~Glasserman and K.-K. Kim.
\newblock Gamma expansion of the {H}eston stochastic volatility model.
\newblock {\em Finance and Stochastics}, 15(2):267--296, 2011.

\bibitem{gm04}
P.~Glasserman and N.~Merener.
\newblock Convergence of a discretization scheme for jump-diffusion processes
  with state-dependent intensities.
\newblock {\em Proc.~Royal Soc.~London~A}, 460:111--127, 2004.

\bibitem{graubner08}
S.~Graubner.
\newblock Multi-level {M}onte {C}arlo {M}ethoden f{\"u}r stochastische
  partielle {D}ifferentialgleichungen.
\newblock Diplomarbeit, TU Darmstadt, 2008.

\bibitem{heinrich98}
S.~Heinrich.
\newblock {M}onte {C}arlo complexity of global solution of integral equations.
\newblock {\em Journal of Complexity}, 14(2):151--175, 1998.

\bibitem{heinrich00}
S.~Heinrich.
\newblock The multilevel method of dependent tests.
\newblock In N.~Balakrishnan, V.B. Melas, and S.~Ermakov, editors, {\em
  Advances in Stochastic Simulation Methods}, pages 47--61. Springer-Verlag,
  2000.

\bibitem{heinrich01}
S.~Heinrich.
\newblock {\em Multilevel {M}onte {C}arlo Methods}, volume 2179 of {\em Lecture
  Notes in Computer Science}, pages 58--67.
\newblock Springer-Verlag, 2001.

\bibitem{heinrich06}
S.~Heinrich.
\newblock {M}onte {C}arlo approximation of weakly singular integral operators.
\newblock {\em Journal of Complexity}, 22(2):192--219, 2006.

\bibitem{hs99}
S.~Heinrich and E.~Sindambiwe.
\newblock {M}onte {C}arlo complexity of parametric integration.
\newblock {\em Journal of Complexity}, 15(3):317--341, 1999.

\bibitem{heston93}
S.I. Heston.
\newblock A closed-form solution for options with stochastic volatility with
  applications to bond and currency options.
\newblock {\em Review of Financial Studies}, 6(2):327--343, 1993.

\bibitem{higham93}
N.J. Higham.
\newblock The accuracy of floating point summation.
\newblock {\em SIAM Journal on Scientific Computing}, 14(4):783--799, 1993.

\bibitem{hsst12}
H.~Hoel, E.~von Schwerin, A.~Szepessy, and R.~Tempone.
\newblock Adaptive multilevel {M}onte {C}arlo simulation.
\newblock In B.~Engquist, O.~Runborg, and Y.-H.R. Tsai, editors, {\em Numerical
  Analysis of Multiscale Computations}, number~82 in Lecture Notes in
  Computational Science and Engineering, pages 217--234. Springer-Verlag, 2012.

\bibitem{kebaier05}
A.~Kebaier.
\newblock Statistical {R}omberg extrapolation: a new variance reduction method
  and applications to options pricing.
\newblock {\em Annals of Applied Probability}, 14(4):2681--2705, 2005.

\bibitem{kp92}
P.E. Kloeden and E.~Platen.
\newblock {\em Numerical Solution of Stochastic Differential Equations}.
\newblock Springer, Berlin, 1992.

\bibitem{marxen10}
H.~Marxen.
\newblock The multilevel {M}onte {C}arlo method used on a {L}{\'e}vy driven
  {SDE}.
\newblock {\em Monte Carlo Methods and Applications}, 16(2):167--190, 2010.

\bibitem{merton76}
R.C. Merton.
\newblock Option pricing when underlying stock returns are discontinuous.
\newblock {\em Journal of Finance}, 3:125--144, 1976.

\bibitem{mss12}
S.~Mishra, C.~Schwab, and J.~Sukys.
\newblock Multi-level {M}onte {C}arlo finite volume methods for nonlinear
  systems of conservation laws in multi-dimensions.
\newblock {\em Journal of Computational Physics}, 231(8):3365--3388, 2012.

\bibitem{mr09}
T.~M{\"u}ller-Gronbach and K.~Ritter.
\newblock Variable subspace sampling and multi-level algorithms.
\newblock In P.~L'Ecuyer and A.~Owen, editors, {\em Monte Carlo and Quasi-Monte
  Carlo Methods 2008}, pages 131--156. Springer-Verlag, 2009.

\bibitem{pb10}
E.~Platen and N.~Bruti-Liberati.
\newblock {\em Numerical Solution of Stochastic Differential Equations With
  Jumps in Finance}.
\newblock Springer, 2010.

\bibitem{schoutens03}
W.~Schoutens.
\newblock {\em L{\'e}vy processes in finance: pricing financial derivatives}.
\newblock Wiley, 2003.

\bibitem{speight09}
A.L. Speight.
\newblock A multilevel approach to control variates.
\newblock {\em Journal of Computational Finance}, 12:1--25, 2009.

\bibitem{speight10}
A.L. Speight.
\newblock Multigrid techniques in economics.
\newblock {\em Operations Research}, 58(4):1057--1078, 2010.

\bibitem{tsgu13}
A.~Teckentrup, R.~Scheichl, M.B. Giles, and E.~Ullmann.
\newblock Further analysis of multilevel {M}onte {C}arlo methods for elliptic
  {PDE}s with random coefficients.
\newblock {\em Numerische Mathematik}, 2013.

\bibitem{xg12}
Y.~Xia and M.B. Giles.
\newblock Multilevel path simulation for jump-diffusion {SDE}s.
\newblock In L.~Plaskota and H.~Wo{\'z}niakowski, editors, {\em Monte Carlo and
  Quasi-Monte Carlo Methods 2010}, pages 695--708. Springer-Verlag, 2012.

\end{thebibliography}

\end{document}